\documentclass[16pt]{article}
\usepackage{amsmath}
\usepackage{graphicx,psfrag,epsf}
\usepackage{enumerate}
\usepackage{natbib}
\usepackage{url} 

\usepackage{comment}

\usepackage{scalerel,stackengine}
\stackMath
\newcommand\reallywidehat[1]{%
\savestack{\tmpbox}{\stretchto{%
  \scaleto{%
    \scalerel*[\widthof{\ensuremath{#1}}]{\kern-.6pt\bigwedge\kern-.6pt}%
    {\rule[-\textheight/2]{1ex}{\textheight}}
  }{\textheight}%
}{0.5ex}}%
\stackon[1pt]{#1}{\tmpbox}%
}

\usepackage[table]{xcolor}
\newcommand{\blind}{0}

\usepackage{soul}
\usepackage{xcolor}     
\usepackage[normalem]{ulem}

\usepackage{calc}
\newsavebox\CBox
\newcommand\hcancel[2][0.5pt]{%
  \ifmmode\sbox\CBox{$#2$}\else\sbox\CBox{#2}\fi%
  \makebox[0pt][l]{\usebox\CBox}%
  \rule[0.5\ht\CBox-#1/2]{\wd\CBox}{#1}}

\begin{document}

\bibliographystyle{agsm}

\def\spacingset#1{\renewcommand{\baselinestretch}%
{#1}\small\normalsize} \spacingset{1}


\if0\blind
{
  \title{\bf Defining a credible interval is not always possible with ``point-null'' priors: A lesser-known correlate of the Jeffreys-Lindley paradox}
  \author{Harlan Campbell and Paul Gustafson \\
  Department of Statistics, University of British Columbia
  }
  \maketitle
} \fi

\if1\blind
{
  \bigskip
  \bigskip
  \bigskip
  \begin{center}
    {\LARGE\bf Title}
\end{center}
  \medskip
} \fi

\bigskip

\abstract{In many common situations, a Bayesian credible interval will be, given the same data, very similar to a frequentist confidence interval, and researchers will interpret these intervals in a similar fashion. However, no predictable similarity exists when credible intervals are based on model-averaged posteriors whenever one of the two {}{}\textcolor{black} {nested} models under consideration is a so called ``point-null''.  Not only can this model-averaged credible interval be quite different than the frequentist confidence interval, in some cases it may be undefined. This is a lesser-known correlate of the Jeffreys-Lindley paradox and is of particular interest given the  popularity of the Bayes factor for testing point-null hypotheses.}

\spacingset{1.45} 

\section{Introduction}

Recently, several Bayesian tests using Bayes factors have been proposed as alternatives to frequentist hypothesis testing; see \citet{heck2020review} for a recent review.  When using the Bayes factor (or the posterior model odds) for testing, it is often recommended that researchers also report parameter estimates and their credible intervals (e.g., \citet{keysers2020using}).  Indeed, following a controversial debate about the strict binary nature of statistical tests, many now call for an additional focus on parameter estimation with appropriate uncertainty estimation; see \citet{wasserstein2016asa}.

\cite{campbell2022bayes} consider how Bayesian testing and estimation can be done in a complimentary manner and conclude that if one reports a Bayes factor comparing two models, then one should also report a \textit{model-averaged} credible interval (i.e., one based on the posterior averaged over the two models under consideration).  Researchers who follow this recommendation can obtain credible intervals congruent with their Bayes factor, thereby obtaining suitable uncertainty estimation.

In many familiar situations, a posterior credible interval will be, given the same data, very similar to a frequentist confidence interval and researchers will interpret these intervals in a similar fashion; see \citet{albers2018credible}.  However, when comparing two models, one of which involves a so-called ``point-null'', it is less clear whether or not such similarity can be assumed.  

Previous work has examined the properties of Bayesian credible intervals and how they relate to frequentist confidence intervals under various prior specifications (e.g., \citet{casella1987reconciling}, {}{}\textcolor{black}{\citet{datta1995priors}}, \citet{greenland2013living}, \citet{held2020bayesian} ).  In this paper, on the basis of a few simple examples, we will examine properties specific to {model-averaged} credible intervals.  We will show that, when one of the two models under consideration is a point-null model, not only can the {}{}\textcolor{black}{model-averaged} credible interval be quite different than the confidence interval, oftentimes, for a desired probability level, it may be undefined.  This is perhaps an unexpected correlate of the Jeffreys-Lindley paradox, the most well known example of the rift between frequentist and Bayesian statistical philosophies; see \citet{wagenmakers2021history}.  The limitations/particularities of working with point-null models are of particular interest given the recent popularity of the Bayes factor for testing point-null hypotheses.

We begin in Section 2 by re-visiting an example of two Normal models considered previously by \citet{wagenmakers2021history} in their discussion of the Jeffreys-Lindley paradox.  In Section 3, we extend this example to consider the consequences of specifying a point-null model.  We conclude in Section 4 with thoughts on the consequences, with respect to parameter estimation, of specifying point-null models.

\section{A mixture of two Normals}


 Let $\theta$ be the parameter of interest for which there are two \textit{a priori} probable models: $M_{0}$ and $M_{1}$, defined by two different priors $\pi_{0}(\theta)$ and $\pi_{1}(\theta)$.   The posterior density which appropriately acknowledges the uncertainty with regards to which of the two models is correct is the mixture density:
\begin{align}
\pi(\theta|data) =& \textrm{Pr}(M_{0}|data)\pi_{0}(\theta|data) +  \textrm{Pr}(M_{1}|data)\pi_{1}(\theta|data),
\label{eq:post_both}
\end{align}
\noindent where the model-specific posteriors, $\pi_{0}(\theta|data)$ and $\pi_{1}(\theta|data)$, are weighted by their posterior model probabilities, $\textrm{Pr}(M_{0}|data)$ and $\textrm{Pr}(M_{1}|data)$; see, for instance, \cite{campbell2022bayes}.  Note that this ``mixture'' posterior is obtained as a result of specifying the ``mixture'' prior:
\begin{equation}
    \pi(\theta) = \textrm{Pr}(M_{0})\pi_{0}(\theta)+\textrm{Pr}(M_{1})\pi_{1}(\theta),
    \label{eq:mixprior}
\end{equation}
where $\textrm{Pr}(M_{0})$ and $\textrm{Pr}(M_{1})$ are the \textit{a priori} model probabilities. 

As an example, consider two \textit{a priori} equally probable Normal models, $M_{0}: \theta \sim N(0, g_{0})$ and $M_{1}: \theta \sim N(0, g_{1})$, such that $\textrm{Pr}(M_{0})=\textrm{Pr}(M_{1})=0.5$.   The prior density functions for the two models are defined as:
\begin{align}
    \pi_{0}(\theta) &= f_{Normal}(\theta, 0, g_{0}), 
    \label{eq:norm_g0}
\end{align}    
\noindent and
\begin{align}
    \pi_{1}(\theta) &= f_{Normal}(\theta, 0, g_{1}), \nonumber
\end{align}
\noindent where $f_{Normal}(x, \mu, \sigma^2)$ is the Normal probability density function evaluated at $x$, with mean parameter $\mu$ and variance parameter $\sigma^2$.  Let $y_{i}$ be the $i$-th data-point, for $i=1,...,n$; let $\bar{y}=\sum_{i=1}^{n}{y_{i}}/n$ be the sample mean; and suppose these data are normally distributed with known unit variance such that:
\begin{align}
    \textrm{Pr}(data|\theta) &= \prod_{i=1}^{n}f_{Normal}(y_{i}, \theta, 1). \nonumber
\end{align}  
Then the Bayes factor is:
\begin{align}
\textrm{BF}_{01} = \sqrt{\frac{1+ng_{1}}{1+ng_{0}}} \times \textrm{exp}\Big(  \frac{(g_{0}-g_{1})nz^{2}}{2(1+ng_{0})(1+ng_{1})}  \Big), \nonumber
\end{align}
where $z=\sqrt{n}\bar{y}$.  The posterior model probabilities can be calculated from the Bayes factor as:
\begin{align}
\textrm{Pr}(M_{0}|data) = \frac{\textrm{Pr}(M_{0})}{ \textrm{Pr}(M_{1})/\textrm{BF}_{01}  + \textrm{Pr}(M_{0})} \quad \textrm{and} \quad \textrm{Pr}(M_{1}|data) = 1- \textrm{Pr}(M_{0}|data).
\label{eq:post_mod_mix}
\end{align}
\noindent Finally, the model specific posteriors are defined as:
\begin{align}
\pi_{j}(\theta|data) =& f_{Normal}\Big(\theta, \frac{zg_{j}}{\sqrt{n}(\frac{1}{n}+g_{j})}, { \frac{g_{j}}{1+g_{j}n}}\Big),  \nonumber
\end{align}
\noindent for $j=0,1$.


Having established all the components of equation (\ref{eq:post_both}), let us now consider how to define a credible interval based on the model-averaged posterior.  An upper one-sided $(1-\alpha)\%$ credible interval is defined as:
\begin{align}
    \textrm{one-sided } (1-\alpha)\%\textrm{CrI} &= [\theta^{*}, \infty), \nonumber
\end{align}
where $\theta^{*}$ satisfies the following equality:
\begin{align}\textrm{Pr}(\theta < \theta^{*}|data) =  \alpha.\label{eq:alpha_thetastar} 
\end{align}
Let us define an equal-tailed two-sided $(1-\alpha)\%$ credible interval from a combination of two upper one-sided intervals as:
\begin{align}
    \textrm{two-sided } (1-\alpha)\%\textrm{CrI} &= [\theta^{l*}, \theta^{u*}), \nonumber
\end{align}
where $\theta^{l*}$  and $\theta^{u*}$ satisfy: $\textrm{Pr}(\theta < \theta^{l*}|data) =  \alpha/2$ and $\textrm{Pr}(\theta < \theta^{u*}|data) =  1-\alpha/2$.  Note that, in our example of two Normal models, these posterior values are calculated as:




%
%
\begin{align}
    \textrm{Pr}(\theta < \theta^{*}|data)     &= \int_{-\infty}^{\theta^{*}}\pi(\theta|data)d\theta  
    = \frac{\int_{-\infty}^{\theta^{*}}\Big(f_{Norm}(    (z-\theta\sqrt{n}),0,1)\times\pi(\theta)\Big)d\theta} {\int_{-\infty}^{\infty}\Big(f_{Norm}((z-\theta\sqrt{n}),0,1) \times \pi(\theta)\Big)d\theta} \nonumber  ,
\end{align}
\noindent where $\pi(\theta)$ is defined as in equation (\ref{eq:mixprior}), {}{}\textcolor{black}{ and the integral in the denominator ensures that the posterior density integrates to one. }

\begin{figure}
      \centering
      \includegraphics[width=15.5cm]{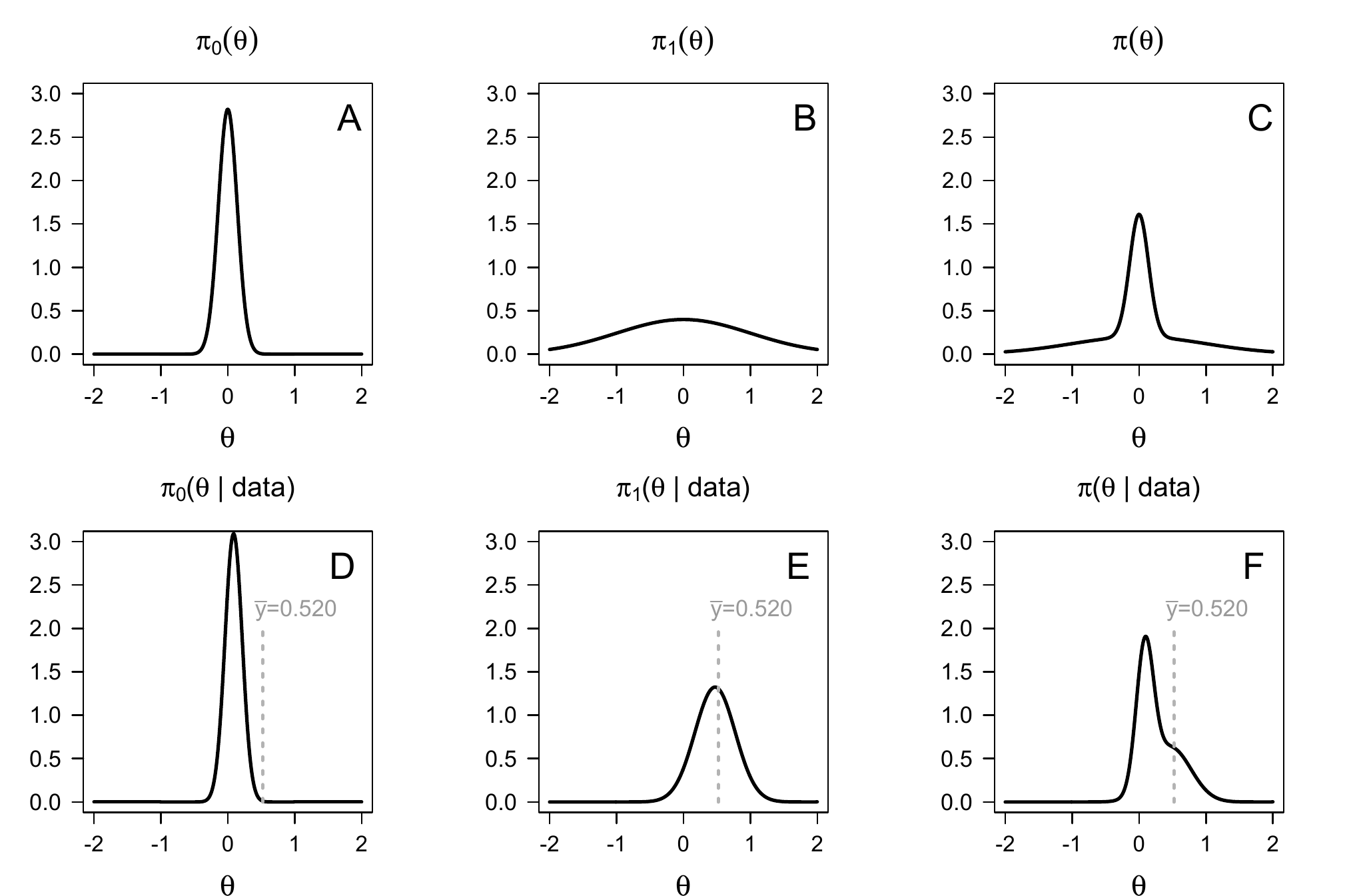}
      \caption{For the ``mixture of two normals'' example ($g_{0}=0.02$ and $g_{1}=1$), panels A, B, and C, plot the $M_{0}$ prior, the $M_{1}$ prior, and the mixture-prior, respectively.  For data with $\bar{y}=0.520$ and $n=10$, panels D, E, and F, plot the $M_{0}$ posterior, the $M_{1}$ posterior, and the model-averaged posterior, respectively.     }
      \label{fig:priorpost1}
  \end{figure}

  {}{}\textcolor{black}{Now suppose} $g_{0}=0.02$, $g_{1}=1$ and that we observe data for which $\bar{y}=1.645/\sqrt{n}$ which corresponds to a $p$-value of $p=0.05$ when using these data to test against the null hypothesis $\textrm{H}_{0}: \theta<0$.   See Figure \ref{fig:priorpost1} which plots priors and posteriors for this scenario with $n=10$.  The lower bound of an upper one-sided $(1-A)$\% confidence interval (CI) will be equal to $\textrm{CI}_{A} = \bar{y} - \frac{Q_{Norm}(1-A)}{\sqrt{n}}$, where $Q_{Norm}()$ is the Normal quantile function.   For instance, for the observed data with $\bar{y}=1.645/\sqrt{n}$, we have $\textrm{CI}_{0.10}= (1.645-1.282)/\sqrt{n}$, such that an upper one-sided 90\% CI will be $ = [0.363/\sqrt{n},\infty)$.  An upper one-sided 95\% CI for these data will be $ [0,\infty)$, since $\textrm{CI}_{0.05} =0$.  How do these frequentist intervals compare to model-averaged Bayesian credible intervals?    {}{}\textcolor{black}{While most literature describing the asymptotic agreement of Bayes and frequentist inferences considers the regime of a fixed true parameter value as $n$ increases, for our purposes it is useful to consider the regime of a fixed $p$-value for a particular point null hypothesis.  Consider two observations.}




  
\begin{figure}
    \centering
    \includegraphics[width=15.5cm]{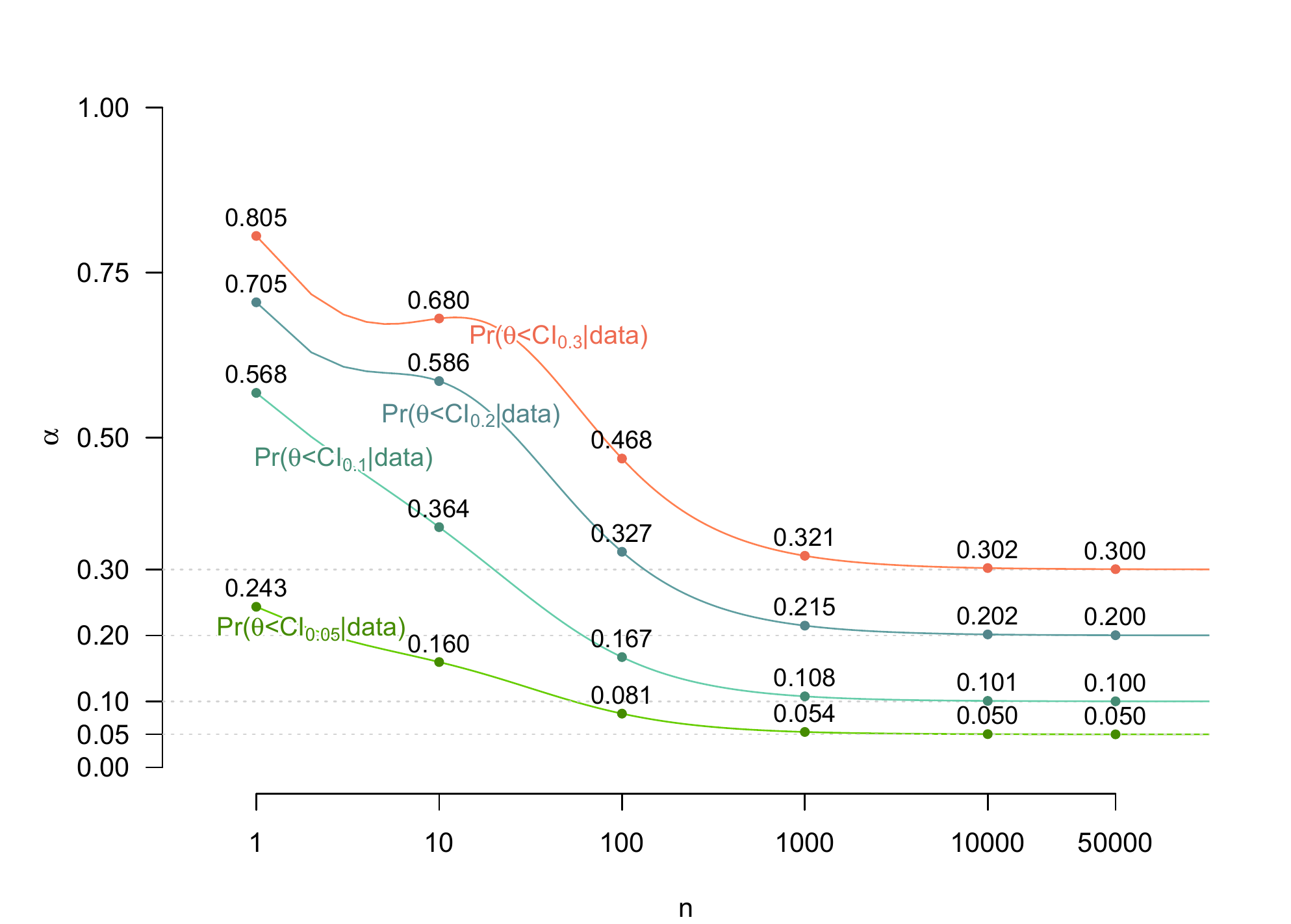}
    \caption{Let $\textrm{CI}_{A}$ be the lower bound of a frequentist upper one-sided (1-A)\% confidence interval.  We consider $\textrm{Pr}(\theta < \textrm{CI}_{A}|data) =  \alpha$ and data corresponding to ($n$, $p$), where $n$ is the sample size and $p$ is the frequentist $p$-value obtained when testing the data against the null hypothesis $\textrm{H}_{0}: \theta<0$. For the normal mixture example with $g_{0}=0.02$ and $g_{1}=1$, and $p=0.05$, we have $\textrm{CI}_{0.05}=0$ and see that, as $n$ increases, $\alpha$ approaches $A$ for $A=0.05,0.10,0.20$ and $0.3$.} 
    \label{fig:alpha_for_p}
\end{figure}

   {}{}\textcolor{black}{First, setting $\theta^{*}=0$ in equation (\ref{eq:alpha_thetastar}), we see that as $n$ increases (and $p=0.05$ remains fixed), the corresponding value of $\alpha$ approaches $p=0.05$: For $n=10$, we obtain $\alpha=0.160$, whereas for $n=10000$, we obtain $\alpha=0.050$;  see how the $\textrm{Pr}(\theta<\textrm{CI}_{0.05}|data)$ curve approaches 0.05 as $n$ increases in Figure \ref{fig:alpha_for_p}.  Second, setting $\theta^{*}= \textrm{CI}_{A}$ in equation (\ref{eq:alpha_thetastar}), we see that as $n$ increases (and $p=0.05$ remains fixed), the corresponding value of $\alpha$  approaches $A$.  In Figure \ref{fig:alpha_for_p}, we plot values of $\alpha$ corresponding to $A=0.05, 0.10, 0.20$, and $0.30$.  One can clearly see that each $\textrm{Pr}(\theta<\textrm{CI}_{A}|data)$ curve tends asymptotically towards $A$.} {}{}\textcolor{black}{One can verify this asymptotic behaviour by re-expressing posterior expectations arising from the specified prior as posterior expectations arising under an improper uniform prior. For completeness, we give the necessary details in the Appendix.}
  

 {}{}\textcolor{black}{Based on the asymptotic behaviour of the posterior in this example, one might reasonably conclude that, with a sufficiently large sample size, the model-averaged credible interval will approximate the frequentist's confidence interval for any $A$ probability level.}  However,  \cite{wagenmakers2021history} argue that, in this scenario,  ``the Jeffreys-Lindley paradox still applies'' indicating that there is indeed a  conflict between Bayesian and frequentist interpretations of the data.
 
\begin{figure}
    \centering
    \includegraphics[width=15.5cm]{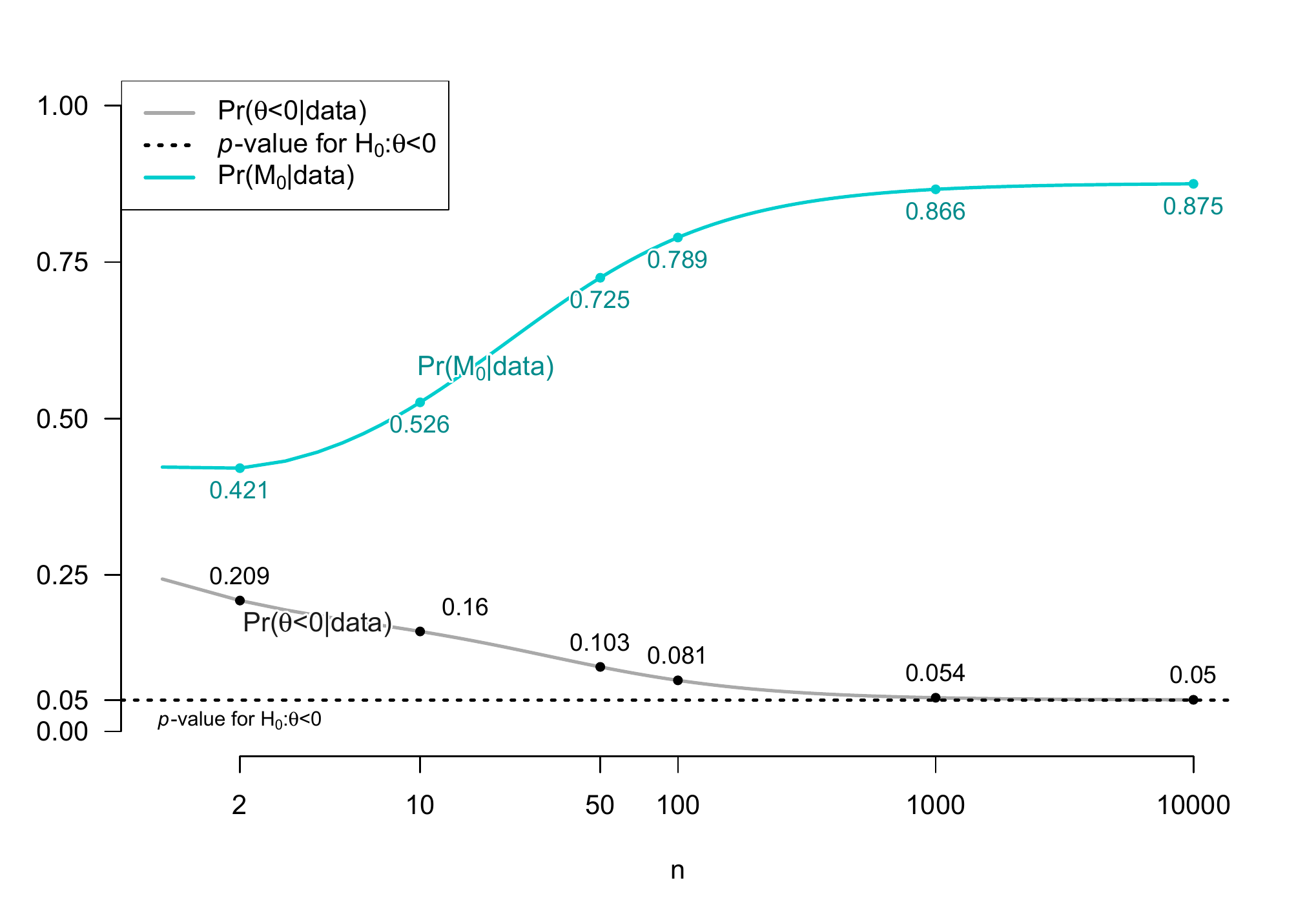}
    \caption{For the normal mixture model example with $g_{0}=0.02$ and $g_{1}=1$, the $\textrm{Pr}(M_{1}|data)$ (blue curve) increases towards 0.876 with increasing $n$, while the value of $\textrm{Pr}(\theta < 0|data)$ (grey line) approaches 0.05 (dashed black line).}
    \label{fig:tost}
\end{figure}

\cite{wagenmakers2021history} explain their reasoning as follows.  From equation (\ref{eq:post_mod_mix}), we calculate $\textrm{lim}_{n\rightarrow\infty}\textrm{Pr}(M_{1}|data) = (1+\sqrt{g_{1}/g_{0}})^{-1} =  (1+1/\sqrt{0.02})^{-1} =0.124$ and $\textrm{lim}_{n\rightarrow\infty}\textrm{Pr}(M_{0}|data)=0.876$.  Therefore, with sufficiently large $n$, we have that $\textrm{Pr}(M_{1}|data) < \textrm{Pr}(M_{0}|data)$ regardless of the data (i.e., regardless of the fixed value of $z=\sqrt{n}\bar{y}$); see Figure \ref{fig:tost}.
 
In this scenario, model selection (i.e., evaluating the relative values of $\textrm{Pr}(M_{0}|data)$ and $\textrm{Pr}(M_{1}|data)$) is not addressing the same question as estimation (i.e., evaluating $\textrm{Pr}(\theta|data)$ to determine which values of $\theta$ are \textit{a posteriori} most likely).  The posterior density of $\theta$ describes one's belief in the probability of different possible values of $\theta$, whereas the posterior model probabilities describe the probability of different data generating processes (DGP) (including the generation of $\theta$).  As such, while it is true that the Jeffreys-Lindley paradox still applies with regards to model selection (i.e., with a sufficiently large sample size and fixed $z$, the Bayesian will inevitably select $M_{0}$), the paradox does not apply when it comes to parameter estimation  (i.e., with a sufficiently large sample size and fixed $z$, the Bayesian will inevitably agree with the frequentist when it comes to estimating $\theta$, with their credible interval approximately equal to the frequentist's confidence interval).  {}{}\textcolor{black}{One way to think about this is to consider the diminishing influence of the prior as the sample size increases and to recall that the confidence interval and the credible interval will agree exactly if one specifies the flat (albeit improper) reference prior, $\pi(\theta) \propto 1$; see the worked examples in \citet{held2020bayesian}.}

In order for the Jeffreys-Lindley paradox to apply to parameter estimation, a point-mass in the prior is required.  We consider this situation in the next Section. 

 

\section{Parameter estimation with a point null}
Consider the same scenario as above but with the null model, $M_{0}$, defined as a so-called ``point-null'' such that the prior density function under $M_{0}$ is:
\begin{align}
    \pi_{0}(\theta) &= \delta_{0}(\theta), 
    \label{eq:pointnullprior}
\end{align}    
\noindent where $\delta_{0}()$ is the Dirac delta function at 0 which can be informally thought of as setting $g_{0}=0$ in equation (\ref{eq:norm_g0}), or alternatively thought of as a probability density function which is zero everywhere except at 0, where it is infinite.  {}{}\textcolor{black}{Note that these are merely informal, intuitive interpretations.} 

  \begin{figure}
      \centering
      \includegraphics[width=15.5cm]{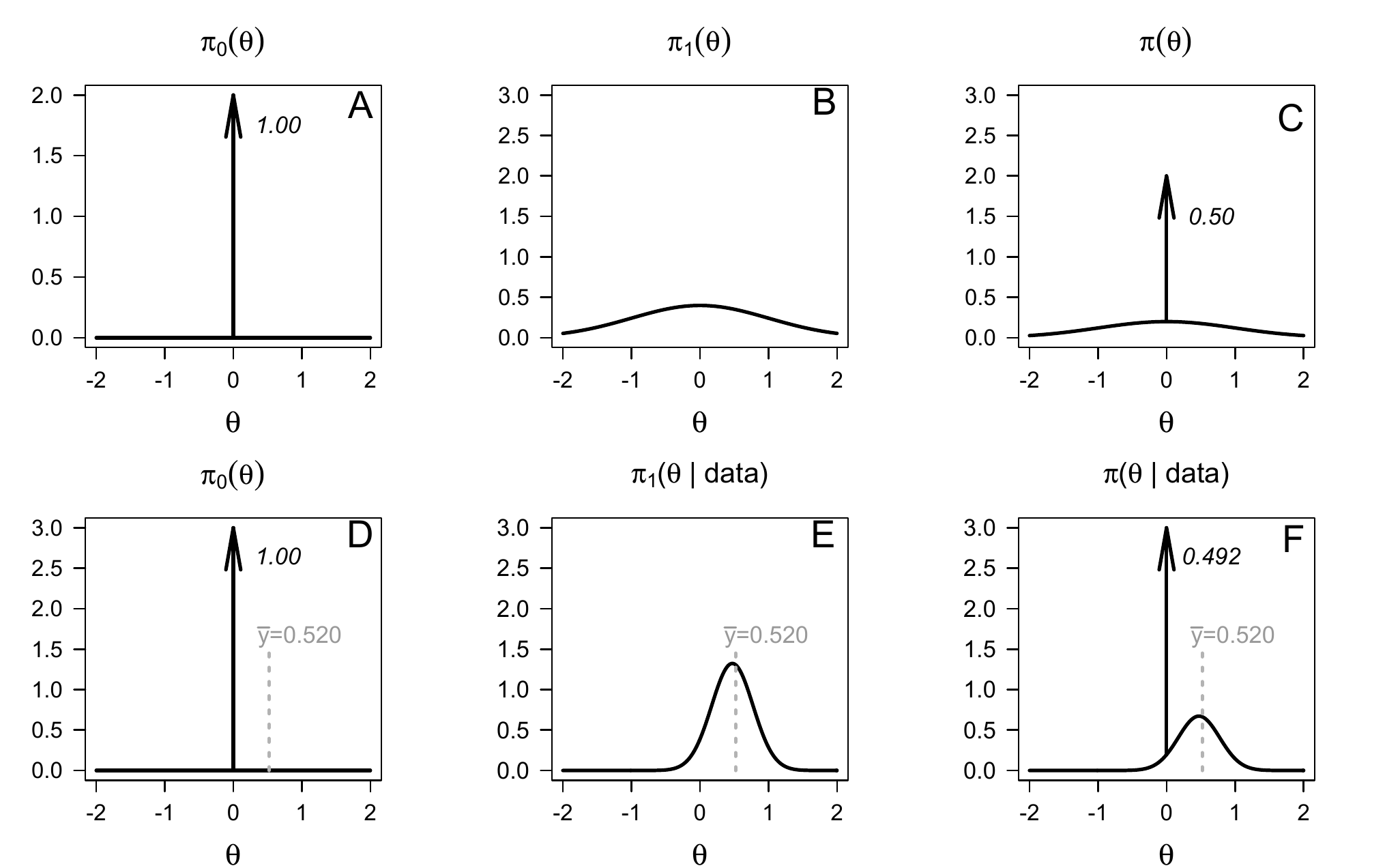}
      \caption{For the ``point-null'' example, panels A, B, and C, plot the $M_{0}$ prior, the $M_{1}$ prior, and the mixture-prior, respectively.  For data with $\bar{y}=0.520$ and $n=10$, panels D, E, and F, plot the $M_{0}$ posterior, the $M_{1}$ posterior, and the model-averaged posterior, respectively.     }
      \label{fig:priorpost2}
  \end{figure}

We now have that $\textrm{Pr}(\theta = 0|data) = \textrm{Pr}(M_{0}|data)$, or equivalently, $\textrm{Pr}(\theta \ne 0|data) = \textrm{Pr}(M_{1}|data)$.  As such, model selection (selecting between $M_{0}$ and $M_{1}$) and null hypothesis testing (selecting between $\textrm{H}_{0}: \theta=0$ and $\textrm{H}_{1}:\theta\ne0$) are equivalent in this scenario.

With the ``point-null'' prior for $M_{0}$ as defined in  (\ref{eq:pointnullprior}), and with $g_{1}=1$, as defined previously, the ``mixture'' prior, $\pi(\theta)$, is recognizable as a ``spike-and-slab'' prior (see \citet{van2021cautionary}) and the Bayes factor is equal to:
\begin{align}
\textrm{BF}_{01} = \sqrt{{1+n}} \times \textrm{exp}\Big(  \frac{-nz^{2}}{2(1+n)}  \Big), \nonumber
\end{align}
The posterior density is nonatomic with a spike (i.e., a discontinuity with infinite density) at 0:
\begin{equation}
    \pi(\theta|data) = \textrm{Pr}(M_{0}|data)\delta_{0}(\theta) + \textrm{Pr}(M_{1}|data)f_{Normal}\Big(\theta, \frac{z}{\sqrt{n}(\frac{1}{n}+1)}, { \frac{1}{1+n}}\Big), \nonumber
\end{equation}
where the posterior model probabilities, $\textrm{Pr}(M_{0}|data)$ and $\textrm{Pr}(M_{1}|data)$, can be calculated from the Bayes factor as in equation (\ref{eq:post_mod_mix}).

\begin{figure}
    \centering
    \includegraphics[width=12.5cm]{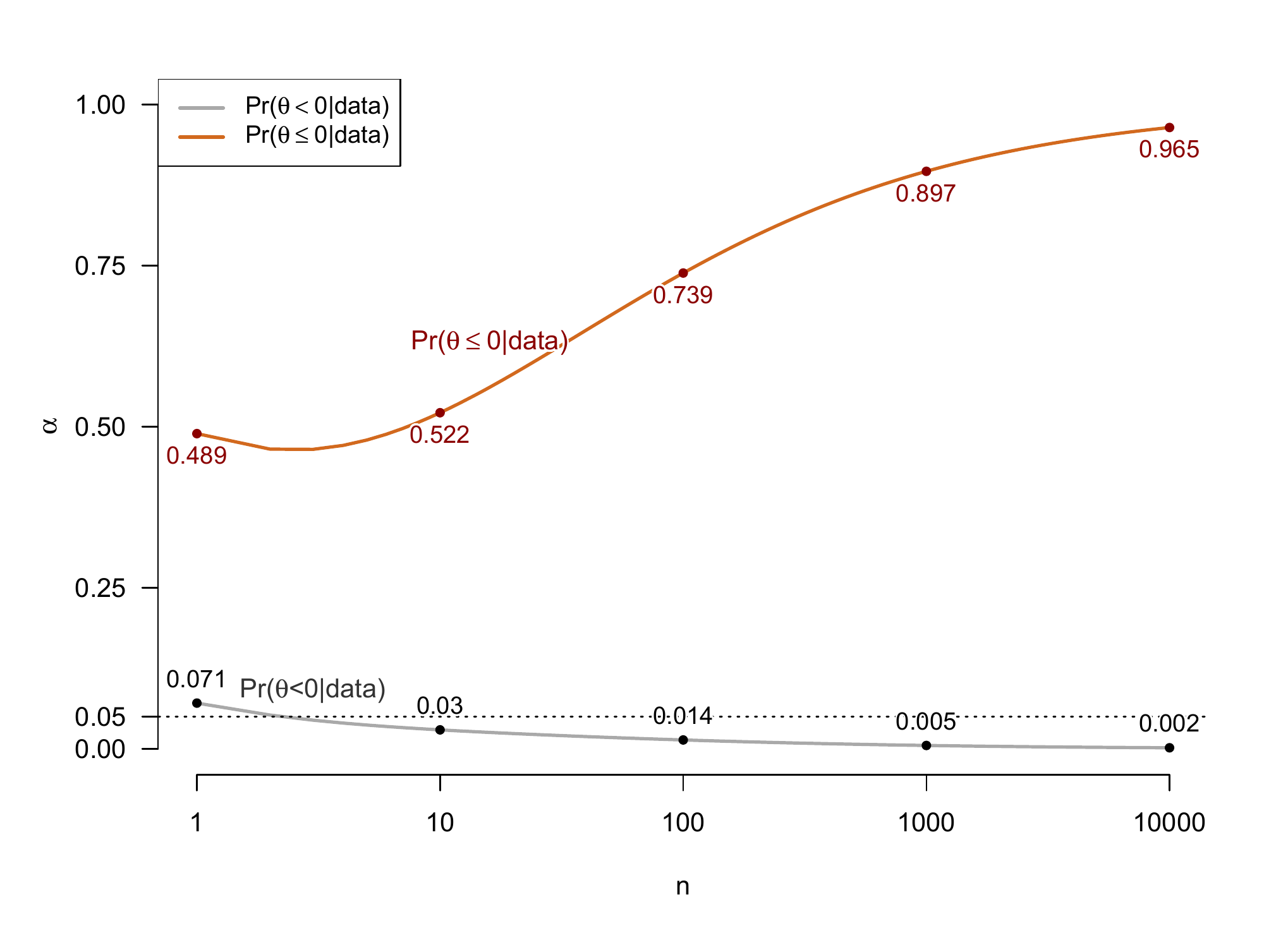}
    \caption{For the hypothetical data with $z=1.645$, as $n$ increases along the horizontal axis, values of $\alpha$ such that $\textrm{Pr}(\theta < 0|data) =  \alpha$ (grey line) and $\textrm{Pr}(\theta \le 0|data) =  \alpha$ (red line) are plotted on the vertical axis.}
    \label{fig:alpha_for_p_point}
\end{figure}

Returning to our hypothetical data with $z=1.645$, we see that for $\theta^{*}=0$, as $n$ increases, $\alpha$ (such that $\textrm{Pr}(\theta < \theta^{*}|data) =  \alpha$) does not approach $p=0.05$ and instead approaches 0: For $n=10$, we obtain $\alpha=0.03$, and for $n=1000$, we obtain $\alpha=0.005$; see trajectory of the grey curve in Figure \ref{fig:alpha_for_p_point}.  Whatsmore, as $n$ increases and $\bar{y}=1.645/\sqrt{n}$ remains fixed, the posterior probability on the ``spike'' at $0$ increases towards infinity such that: $\textrm{lim}_{n\rightarrow \infty}\textrm{Pr}(M_{0}|data) = 1$; as famously {}{}\textcolor{black}{emphasized} by \citet{lindley1957statistical} and originally demonstrated by \citet{jeffreys1935some}.

Perhaps even more puzzling is that, for fixed $\alpha=0.05$, there is simply no corresponding value of $\theta^{*}$ (such that $\alpha = \textrm{Pr}(\theta < \theta^{*}|data)$) for any $n>2$.   For $n=2$ we can define $\theta^{*}=-0.0163$, such that  $\textrm{Pr}(\theta < -0.0163|data)=0.05$.  However, for $n=3$, a precise value of $\theta^{*}$ cannot be defined since, due to the discontinuity in the posterior, we have: $\textrm{Pr}(\theta < 0|data) =0.045<\alpha$, and $\textrm{Pr}(\theta \le 0|data) =0.465>\alpha$.  For $n=10$ the gap is even wider: $\textrm{Pr}(\theta < 0|data) =0.030<\alpha$ and $\textrm{Pr}(\theta \le 0|data) =0.522>\alpha$.  Figure \ref{fig:alpha_for_p_point} plots these numbers for increasing values of $n$.  As a consequence, it is no longer the case that, with a sufficiently large sample size, a Bayesian's credible interval will approximate a frequentist's confidence interval.  In fact, for certain values of $\alpha$ and $n$, calculating a credible interval is not even possible. 



\begin{figure}
    \centering
    \includegraphics[width=12.5cm]{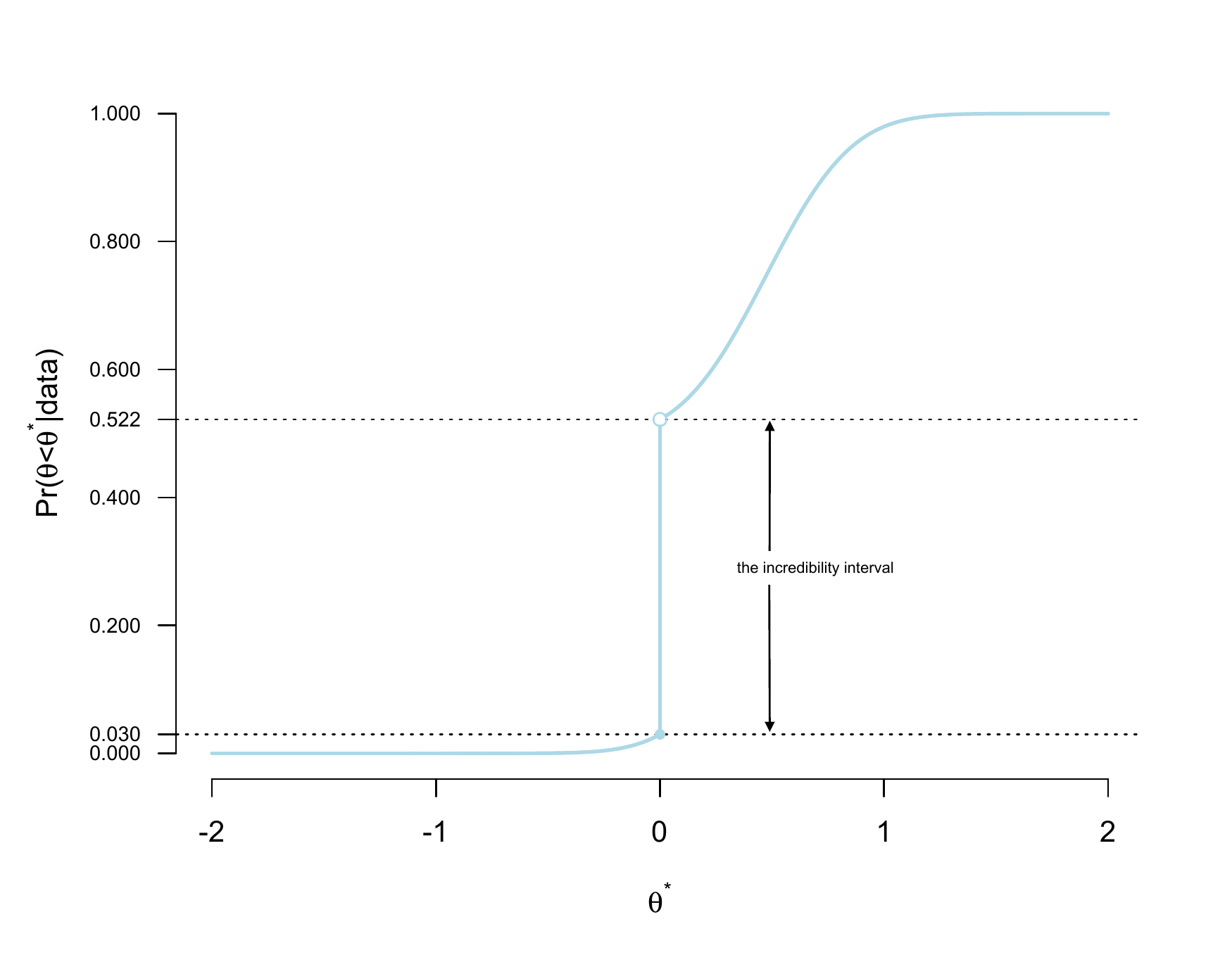}
    \caption{For the hypothetical data with $z=1.645$ and $n=10$, the plotted line corresponds to the cumulative distribution function of the posterior (i.e., $\textrm{Pr}(\theta < \theta^{*}|data)$) for increasing  values of $\theta^{*}$.}
    \label{fig:cdfplot}
\end{figure}

\begin{figure}
    \centering
    \includegraphics[width=15.5cm]{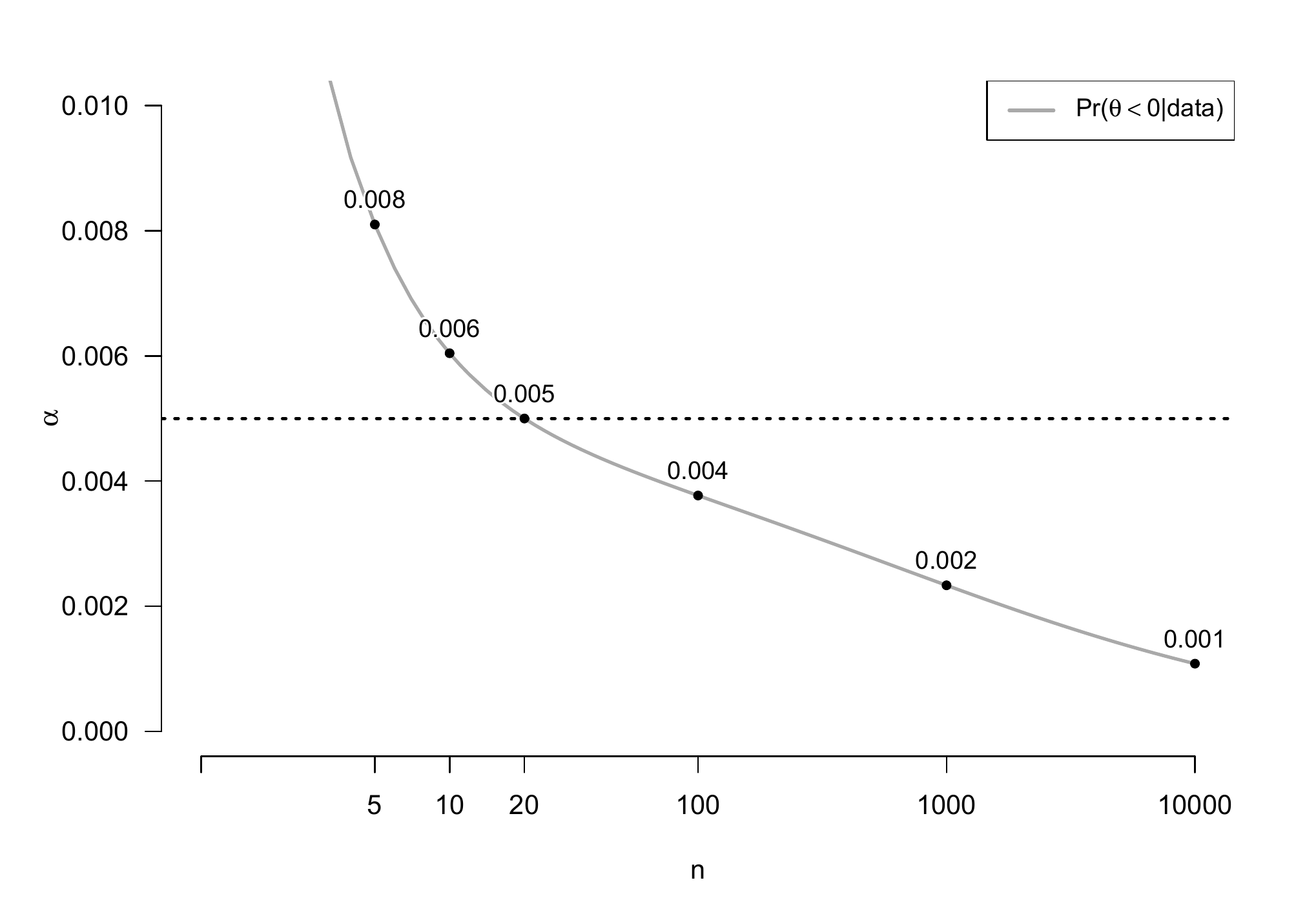}
    \caption{With data where $\bar{y}=2.575/\sqrt{n}$, as $n$ increases, the lower bound of the incredibility interval (the solid line) decreases towards zero.  As a consequence, determining a value of $\theta^{*}$ such that $    \textrm{Pr}(\theta < \theta^{*}|data) =  \alpha$, when $\alpha=0.005$ (the dotted line) is only possible for $n<20$.}
    \label{fig:alpha_0005}
\end{figure}

In general, determining a specific value of $\theta^{*}$ for a given value of $\alpha$ (such that $\alpha = \textrm{Pr}(\theta < \theta^{*}|data)$) is only possible for values of $\alpha$ outside of the ``incredibility interval'':
\begin{equation}
\Big[ \Big(\textrm{Pr}({\theta}<0|data, M_{1}) \textrm{Pr}(M_{1}|data)\Big), \Big(\textrm{Pr}({\theta}<0|data, M_{1}) \textrm{Pr}(M_{1}|data) + \textrm{Pr}(M_{0}|data)\Big) \Big]. \nonumber
\end{equation}
  {}{}\textcolor{black}{The bounds of the ``incredibility interval'' are the limits of the ``jump'' in the cumulative distribution function of the posterior, i.e., the values ranging between $\textrm{Pr}(\theta<0|\textrm{data})$ and $\textrm{Pr}(\theta\le0|\textrm{data})$.  In Figure \ref{fig:cdfplot}, we plot the cumulative distribution function of the posterior  for  hypothetical data with $z=1.645$ and $n=10$.  In this situation, the ``incredibility interval'' equals $[\textrm{Pr}(\theta < 0|data), \textrm{Pr}(\theta \le 0|data)] = [0.03, 0.522]$.} In Figure \ref{fig:alpha_for_p_point}, the lower grey curve corresponds to the lower bound of the incredibility interval and the upper red curve corresponds to the upper bound.  Notably, since $\textrm{lim}_{n\rightarrow \infty}\textrm{Pr}(M_{0}|data) = 1$ and $\textrm{lim}_{n\rightarrow \infty}\textrm{Pr}(M_{1}|data) = 0$, the width of the incredibility interval increases as $n$ increases.  As a result, determining a precisely $\alpha$-level value of $\theta^{*}$  such that $\alpha = \textrm{Pr}(\theta < \theta^{*}|data)$, becomes increasingly impossible as $n$ grows large.  This is true regardless of the data; see Figure \ref{fig:alpha_0005} for values of the lower bound obtained with data where $\bar{y}=2.575/\sqrt{n}$ (data for which one obtains a $p$-value of $p=0.005$ when testing against $\textrm{H}_{0}: \theta<\theta_{0}$).

When $\alpha$ is inside the incredibility interval, there remains an unconventional way for defining a $(1-\alpha)\%$ credible interval. In order to establish a correct value for $\theta^{*}$ such that $\textrm{Pr}(\theta < \theta^{*}|data) =  \alpha$ (over repeated samples) one defines $\theta^{*}$ stochastically such that 
\begin{align}
 \theta^{*}=
\begin{cases}
    0, &\text{with probability } \gamma; \textrm{ and} \\
    0+\epsilon, &\text{with probability } 1-\gamma, 
\end{cases}
\label{eq:sto1}
\end{align}
\noindent where:
\begin{align}\gamma=\frac{\alpha-\textrm{Pr}(\theta\le0|data)}{ \textrm{Pr}(\theta<0|data)-\textrm{Pr}(\theta\le0|data)}, \nonumber
\end{align}
and $\epsilon$ is an arbitrarily small number.

Returning to our example data with $\bar{y}=1.645/\sqrt{n}$, we note that, for $n=10$, $\textrm{Pr}(\theta < 0|data)= 0.030$ and $\textrm{Pr}(\theta\le0|data)= 0.522$.  As such, for $\alpha=0.05$ (which is inside the incredibility interval of [0.030, 0.522]), we define $\theta^{*}$ as: 
\begin{align}
 \theta^{*}=
\begin{cases}
    0, &\text{with probability } \gamma= 0.959; \textrm{ and} \\
    0+\epsilon, &\text{with probability } (1-\gamma)=0.041. \nonumber
\end{cases}
\end{align}
{}{}\textcolor{black}{Defining $\theta^{*}$ in this way will guarantee that $\textrm{Pr}(\theta < \theta^{*}|data) =  0.05$.  One way to  think about this is to consider the various values of $\theta$ that, over a researcher's lifetime give rise to the various datasets they analyse.  Across all of these studies, the average posterior probability content of the $[\theta^{*}, \infty)$ interval will be 0.95.  Moreover, if the model-averaged prior does in fact correspond to the true data generating mechanism, we can be assured that, amongst all of the researcher's studies for which $z=1.645$, 95\% of these were the result of a $\theta$ value from inside of their interval. Furthermore, since this is true for any arbitrary value of $z$ and any arbitrary value of $\alpha$, then we have that $\textrm{Pr}\Big(\theta_{j} \in [\theta^{*}, \infty)\Big|z_{j}\Big)=1-\alpha$, where $\theta_{j}$ and $z_{j}$ are values obtained from a joint draw from the amalgamation of the prior and statistical model (i.e., the data generating mechanism).}


As another example, suppose $n=100$ and $\bar{y}=2.054/\sqrt{n} = 0.2054$ which corresponds to a $p$-value of $p=0.04$ when using the data to test against the null hypothesis $\textrm{H}_{0}: \theta=0$, and a $p$-value of $p=0.02$ when using the data to test against the null hypothesis $\textrm{H}_{0}: \theta<0$.  One can easily calculate an upper one-sided frequentist 95\% confidence interval for these data equal to: $[\bar{y} - 1.645/\sqrt{n}, \infty) = [0.040, \infty)$, which clearly excludes 0.  However, one cannot calculate an upper one-sided 95\% credible interval since $\alpha=0.05$ is within the incredibility interval for this data: [0.009, 0.564].  The closest one can do is to calculate an upper one-sided 99.1\% credible equal to: $[0, \infty)$ which includes 0, or calculate an upper one-sided 43.6\% credible interval equal to $(0, \infty)$ which excludes 0.   The only way to define an upper one-sided interval with exactly 95\% probability of including the true value of $\theta$ (over repeated samples) is to do so stochastically as equal to: $[\theta^{*}, \infty)$, where $\theta^{*}=0$ with probability $\gamma = (0.050-0.564)/(0.009-0.564) = 0.926$, and $\theta^{*}=0+\epsilon$ with probability $1-\gamma = 0.074$.  
  
%

{}{}\textcolor{black}{We are not seriously suggesting that researchers define credible intervals in this bizarre stochastic way.  We simply wish to demonstrate that this is the only way one can correctly define the credible interval from a posterior with point masses.  When model-averaged posteriors involve point-null models, credible intervals must therefore be approached and interpreted with the utmost caution.}  The issue only gets thornier as the sample size increases.


For a very very large $n$ it is possible that both $\alpha/2$ and $(1-\alpha/2)$ are within the incredibility interval.  In this case, the equal-tailed two-sided $(1-\alpha)\%$ credible interval must be defined in an even more bizarre way.  When both $\alpha/2$ and $(1-\alpha/2)$ are both in the incredibility interval, the credible interval must be defined stochastically as either a single point or as an entirely empty interval:
\begin{align}
 (1-\alpha)\%\textrm{CrI}=
\begin{cases}
    [0], &\text{with probability } \psi; \textrm{ and} \\
    	\emptyset , &\text{with probability } (1-\psi), 
\end{cases}
\label{eq:sto2}
\end{align}
where 
\begin{align}
\psi = \frac{\textrm{Pr}(\theta=0|data)-\alpha}{2\times \textrm{Pr}(\theta=0|data)-1}. \nonumber
\end{align}

\begin{figure}
     \centering
     \includegraphics[width=15cm]{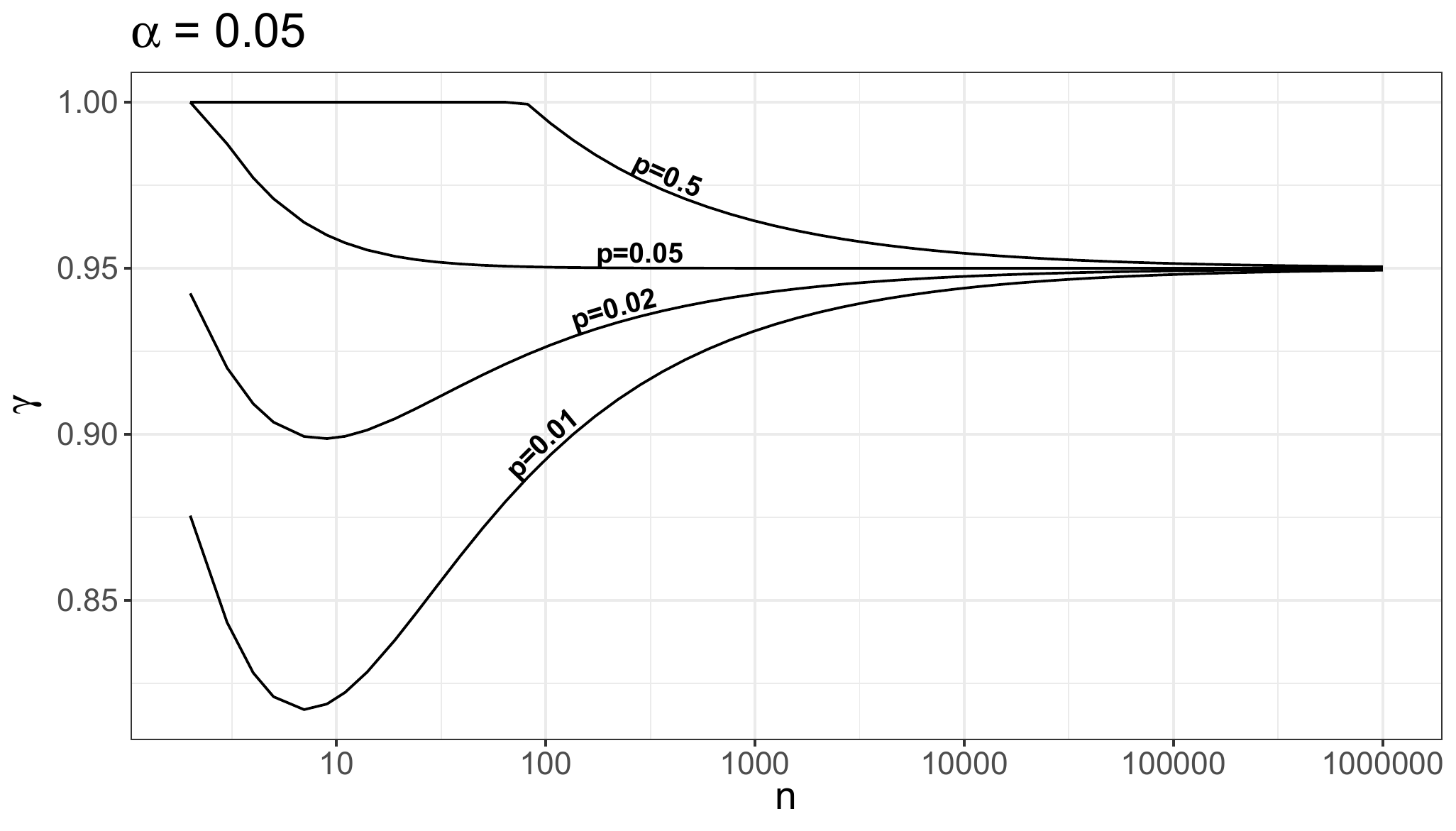}
     \caption{ Each line corresponds to observing data corresponding to a $p$-value of $p$ when testing against $\textrm{H}_{0}:\theta<0$.}
     \label{fig:gamma_plot}
 \end{figure}

 {}{}\textcolor{black}{To be clear, the ``stochastic credible interval'' is not defined in equations (\ref{eq:sto1}) and (\ref{eq:sto2}) to ensure that it has a certain (asymptotic) coverage.  Rather it is defined in the only possible way such that (over repeated samples) the boundaries of the interval contain the correct amount of posterior mass (as  required by the definition in equation (\ref{eq:alpha_thetastar})).  
 As such, it may not be immediately obvious that, when we look at the asymptotic behaviour of these stochastic credible intervals, we see that the Jeffereys-Lindley paradox reduces the data to be entirely inconsequential {}{}\textcolor{black}{(at least when assuming a fixed $p$-value)}.  Indeed,}  as $n$ increases,  both $\gamma$ and $\psi$  approach $1-\alpha$ since:
 \begin{align}
 \textrm{lim}_{n \rightarrow \infty}\gamma &=  \textrm{lim}_{n \rightarrow \infty}\Big(\frac{\alpha-\textrm{Pr}(\theta\le\theta_{0}|data)}{ \textrm{Pr}(\theta<\theta_{0}|data)-\textrm{Pr}(\theta\le\theta_{0}|data)}\Big) \nonumber\\
  &= \Big(\frac{\alpha-1}{-1}\Big) \nonumber \\
 &=  1-\alpha, \nonumber
 \end{align}
 and:
  \begin{align}
 \textrm{lim}_{n \rightarrow \infty}\psi &=  \textrm{lim}_{n \rightarrow \infty}\Big(\frac{\textrm{Pr}(\theta=0|data)-\alpha}{2\times \textrm{Pr}(\theta=0|data)-1}\Big) \nonumber\\
  &= \Big(\frac{1-\alpha}{2-1}\Big) \nonumber \\
 &=  1-\alpha. \nonumber
 \end{align}
 Therefore, for sufficiently large $n$ and $z$ remaining constant, the probability that one will  exclude 0 from a $(1-\alpha)$\%credible interval will equal $\alpha$ regardless of the data; see Figure \ref{fig:gamma_plot}.  While this may strike one as paradoxical, it is entirely congruent with the wildly-known consequence of the Jeffereys-Lindley paradox: As $n$ increases and $z$ is fixed, the probability of selecting $M_{0}$ will go to 1.

 \section{Conclusion}
 
We demonstrated that when one of the two models under consideration is a point-null model, not only can {}{}\textcolor{black}{ a model-averaged} credible interval be rather different than the frequentist confidence interval, oftentimes it will be simply undefined (at least in a conventional sense). {}{}\textcolor{black}{As a consequence, it may be tempting to compare (e.g., using the Bayes factor) two \textit{a priori} probable models,  $M_{0}$ and $M_{1}$,  for the purpose of model selection, but then simply report the uncertainty about $\theta$, conditional on  $M_{1}$ being unquestionably true.  We caution that this strategy, while seemingly straightforward, will lead to unavoidable inconsistencies between one's priors and posteriors.  \citet{campbell2022bayes} explain in detail why disregarding  $M_{0}$ ``for the purpose of parameter estimation''  \citep{wagenmakers2020overwhelming} is inadvisable; see also \citet{tendeiro2019review}. }

{}{}\textcolor{black}{Some researchers may be happy to avoid model selection entirely and may see no reason to entertain point-null priors (e.g., \citet{gelman1995avoiding}: ``realistic prior distributions in social science do not have a mass of probability at zero'' [...]  ``we believe model selection to be relatively unimportant compared to the task of constructing realistic models that agree with both theory and data.'').}  However, if researchers truly believe that there is a non-zero prior probability that the parameter of interest is precisely zero (and this prior probability is equal to the value assigned to $\textrm{Pr}(M_{0})$), Bayesian testing  with a point-null will be optimal in the sense of minimizing the expected loss (with respect to a joint distribution of the data and parameters); see \citet{berger2013statistical}.  {}{}\textcolor{black}{These researchers should be aware that, while perhaps optimal, Bayesian testing with a point-null can lead to rather unexpected asymptotic behaviour.  There will still be credible intervals; it is just that, as a consequence of the discontinuity in the model averaged posterior, certain specific credible intervals do not exist.  Some researchers might therefore wish to explore alternative means of conveying the uncertainty surrounding the parameter of interest (e.g., \citet{wagenmakers2020support}, \citet{rice2022expressing}). } 

{}{}\textcolor{black}{One limitation of this work is that we only considered univariate models where one wishes to define a credible interval for a single parameter of interest.  However, the ideas we discussed also apply to multivariate settings where one wishes to define credible sets and where there may be several different nested models under consideration.  For instance, researchers using Bayes factors in multiple regression models \citep{rouder2012default} should be aware that it may be impossible to define certain model-averaged credible intervals/sets for the regression coefficients. A second limitation is that we did not consider how the undefinability of specific credible intervals will also occur in discrete parameter models.  In such cases, specific confidence intervals will also be undefined \citep{tingley1993note, berger2013statistical}, so while they may both be puzzled, Bayesians and frequentists should at least agree in their inability to define an uncertainty interval!}

Finally, we note that the consequences of the Jeffereys-Lindley paradox on model selection (and null hypothesis testing) are often understood as ``intuitive'' and not necessarily unfavourable: When sample sizes are very large, researchers might indeed prefer to sacrifice some power in order to lower the probability of a type I error, a trade-off that occurs necessarily when testing a point-null hypothesis with the Bayes factor; see \citet{pericchi2016adaptative} and \citet{wagenmakers2021history}.  {}{}\textcolor{black}{Indeed, the benefits of such a trade-off are routinely discussed by frequentists and Bayesians alike (e.g., \citet{leamer1978specification}: ``from every reasonable viewpoint the significance level should be a decreasing function of sample size''; and recently,  \citet{wulff2023and}: ``From a Neyman-Pearson perspective, it is logical that $\alpha$ should be a decreasing function of the sample size.'').}  However, the consequences of the Jeffereys-Lindley paradox on parameter estimation --specifically with regards to model-averaged credible intervals and the inability to define these for certain probability levels-- were previously less well understood, and certainly strike us as less intuitive.

\bibliography{truthinscience}

@article{tingley1993note,
  title={A note on obtaining confidence intervals for discrete parameters},
  author={Tingley, Maureen and Li, Chunhui},
  journal={The {A}merican {S}tatistician},
  volume={47},
  number={1},
  pages={20--23},
  year={1993},
  publisher={Taylor \& Francis}
}

@article{wagenmakers2020support,
  title={The support interval},
  author={Wagenmakers, Eric-Jan and Gronau, Quentin F and Dablander, Fabian and Etz, Alexander},
  journal={Erkenn},
  volume={87},
  pages={589--601},
  year={2022},
  publisher={Springer}
}

@book{berger2013statistical,
  title={Statistical decision theory and {B}ayesian analysis},
  author={Berger, James O},
  year={1985},
  publisher={Springer Science \& Business Media}
}

@article{lindley1957statistical,
  title={A statistical paradox},
  author={Lindley, Dennis V},
  journal={Biometrika},
  volume={44},
  number={1/2},
  pages={187--192},
  year={1957},
  publisher={JSTOR}
}

@article{greenland2013living,
  title={Living with $p$-values: {R}esurrecting a {B}ayesian Perspective on Frequentist Statistics},
  author={Greenland, Sander and Poole, Charles},
  journal={Epidemiology},
  pages={62--68},
  year={2013},
  publisher={JSTOR}
}

@article{casella1987reconciling,
  title={Reconciling {B}ayesian and frequentist evidence in the one-sided testing problem},
  author={Casella, George and Berger, Roger L},
  journal={Journal of the American Statistical Association},
  volume={82},
  number={397},
  pages={106--111},
  year={1987},
  publisher={Taylor \& Francis}
}

@article{campbell2022bayes,
  title={Bayes factors and posterior estimation: Two sides of the very same coin},
  author={Campbell, Harlan and Gustafson, Paul},
  journal={arXiv preprint arXiv:2204.06054},
  year={2022}
}

@article{wagenmakers2021history,
  title={History and nature of the {J}effreys-{L}indley paradox},
  author={Wagenmakers, Eric-Jan and Ly, Alexander},
  journal={arXiv preprint arXiv:2111.10191},
  year={2021}
}

@misc{wagenmakers2020overwhelming,
  title={Overwhelming Evidence for Vaccine Efficacy in the {P}fizer Trial: An Interim {B}ayesian Analysis},
  author={Wagenmakers, Eric-Jan and Gronau, Quentin Frederik},
  year={2020},
  howpublished={PsyArXiv}
}

@article{pericchi2016adaptative,
  title={Adaptative significance levels using optimal decision rules: balancing by weighting the error probabilities},
  author={Pericchi, Luis and Pereira, Carlos},
  journal={Brazilian Journal of Probability and Statistics},
  volume={30},
  number={1},
  pages={70--90},
  year={2016},
  publisher={Brazilian Statistical Association}
}

@article{keysers2020using,
  title={Using {B}ayes factor hypothesis testing in neuroscience to establish evidence of absence},
  author={Keysers, Christian and Gazzola, Valeria and Wagenmakers, Eric-Jan},
  journal={Nature {N}euroscience},
  volume={23},
  number={7},
  pages={788--799},
  year={2020},
  publisher={Nature Publishing Group}
}

@article{albers2018credible,
  title={Credible confidence: {A} pragmatic view on the frequentist vs {B}ayesian debate},
  author={Albers, Casper J and Kiers, Henk AL and van Ravenzwaaij, Don},
  journal={Collabra: {P}sychology},
  volume={4},
  number={1},
  year={2018},
  publisher={University of California Press}
}

@article{wasserstein2016asa,
  title={The {A}{S}{A} statement on $p$-values: context, process, and purpose},
  author={Wasserstein, Ronald L and Lazar, Nicole A},
  journal={The {A}merican {S}tatistician},
  volume={70},
  number={2},
  pages={129--133},
  year={2016},
  publisher={Taylor \& Francis}
}

@article{tendeiro2019review,
  title={A review of issues about null hypothesis {B}ayesian testing.},
  author={Tendeiro, Jorge N and Kiers, Henk AL},
  journal={Psychological {M}ethods},
  volume={24},
  number={6},
  pages={774},
  year={2019},
  publisher={American Psychological Association}
}

@article{rice2022expressing,
  title={Expressing regret: a unified view of credible intervals},
  author={Rice, Kenneth and Ye, Lingbo},
  journal={The {A}merican {S}tatistician},
  volume={76},
  number={3},
  pages={248--256},
  year={2022},
  publisher={Taylor \& Francis}
}

@article{rouder2012default,
  title={Default {B}ayes factors for model selection in regression},
  author={Rouder, Jeffrey N and Morey, Richard D},
  journal={Multivariate Behavioral Research},
  volume={47},
  number={6},
  pages={877--903},
  year={2012},
  publisher={Taylor \& Francis}
}

@article{gelman1995avoiding,
  title={Avoiding model selection in {B}ayesian social research},
  author={Gelman, Andrew and Rubin, Donald B},
  journal={Sociological methodology},
  volume={25},
  pages={165--173},
  year={1995},
  publisher={JSTOR}
}

@article{van2021cautionary,
  title={A cautionary note on estimating effect size},
  author={van den Bergh, Don and Haaf, Julia M and Ly, Alexander and Rouder, Jeffrey N and Wagenmakers, Eric-Jan},
  journal={Advances in Methods and Practices in Psychological Science},
  volume={4},
  number={1},
  pages={2515245921992035},
  year={2021},
  publisher={Sage Publications Sage CA: Los Angeles, CA}
}

@article{datta1995priors,
  title={On priors providing frequentist validity for {B}ayesian inference},
  author={Datta, Gauri Sankar and Ghosh, Jayanta Kumar},
  journal={Biometrika},
  volume={82},
  number={1},
  pages={37--45},
  year={1995},
  publisher={Oxford University Press}
}

@article{wulff2023and,
  title={How and why alpha should depend on sample size: A {B}ayesian-frequentist compromise for significance testing},
  author={Wulff, Jesper N and Taylor, Luke},
  year={2023}
}

@book{leamer1978specification,
  title={Specification searches: {A}d hoc inference with nonexperimental data},
  author={Leamer, Edward E},
  volume={53},
  year={1978},
  publisher={John Wiley \& Sons Incorporated}
}

@inproceedings{jeffreys1935some,
  title={Some tests of significance, treated by the theory of probability},
  author={Jeffreys, Harold},
  booktitle={Mathematical proceedings of the Cambridge philosophical society},
  volume={31},
  number={2},
  pages={203--222},
  year={1935},
  organization={Cambridge University Press}
}

@incollection{held2020bayesian,
  title={Bayesian Tail Probabilities for Decision Making},
  author={Held, Leonhard},
  booktitle={Bayesian Methods in Pharmaceutical Research},
  pages={53--73},
  year={2020},
  publisher={CRC Press Taylor \& Francis Group}
}

@article{heck2020review,
  title={A review of applications of the {B}ayes factor in psychological research},
  author={Heck, Daniel W and Boehm, Udo and B{\"o}ing-Messing, Florian and B{\"u}rkner, Paul-Christian and Derks, Koen and Dienes, Zoltan and Fu, Qianrao and Gu, Xin and Karimova, Diana and Kiers, Henk AL and others},
  journal={Psychological Methods},
  year={2022},
  publisher={American Psychological Association}
}

\pagebreak

\section{Appendix}
  {}{}\textcolor{black}{
To verify the behavior seen in Figure 2, 
consider the asymptotic regime with $\bar{y} = a + b n^{-1/2}$, for some $a$, 
and some $b>0$.
This corresponds to a fixed $p$-value against $\textrm{H}_{0}: \theta < a$.
We are interested in the posterior probability content of the CI $(\bar{y}-kn^{-1/2},\infty)$ $=$
$(a+(b-k)n^{-1/2},\infty)$,
with $k>0$ chosen to give the desired coverage, i.e., a coverage probability of $\Phi(k)$, where $\Phi()$ is the standard normal CDF.
We presume the prior density $\pi()$ on $\theta$ to be continuous, with $\pi(a)>0$.
In what follows we let 
Let $\textrm{E}^{*}_{n}$ indicate expectation with respect to $\theta \sim N(a+bn^{-1/2},n^{-1})$ (so the posterior on $\theta$ that would arise under a locally uniform prior).   And let $\textrm{E}^{*}$ indicate expectation with respect to $S \sim N(0,1)$ (as one could get by standardizing the distribution above).  
}
  {}{}\textcolor{black}{
With this set-up, the posterior probability of interest is the complement of:
\begin{eqnarray*}
\textrm{Pr} \left( \left. \theta < a + \frac{b-k}{n^{1/2}} \right| \bar{Y}= a + \frac{b}{n^{1/2}} \right)
&=& \int_{-\infty}^{a + (b-k)n^{-1/2}}\pi(\theta| \bar{Y}= a + {b}{n^{-1/2}})d\theta\\
           &=& \frac{\int_{-\infty}^{a + (b-k)n^{-1/2}}\Big(f_{Norm}( a+bn^{-1/2},\theta,n^{-1})\times\pi(\theta)\Big)d\theta} {\int_{-\infty}^{\infty}\Big(f_{Norm}(a+bn^{-1/2},\theta,n^{-1}) \times \pi(\theta)\Big)d\theta}\\
      &=& \frac{\int_{-\infty}^{a + (b-k)n^{-1/2}}\Big(f_{Norm}( \theta, a+bn^{-1/2},n^{-1})\times\pi(\theta)\Big)d\theta} {\int_{-\infty}^{\infty}\Big(f_{Norm}(\theta, a+bn^{-1/2},n^{-1}) \times \pi(\theta)\Big)d\theta}\\
&=& 
\frac{\textrm{E}^{*}_{n} \left\{ I_{(-\infty, a + (b-k)n^{-1/2})}(\theta) \pi(\theta) \right\} }
{\textrm{E}^{*}_{n}\left\{ \pi(\theta) \right\}} \\
& = &
\frac{\textrm{E}^{*} \left\{I_{(-\infty, -k)}(S) \pi\left(a+(b+S)n^{-1/2}\right) \right\}}
{\textrm{E}^{*} \left\{ \pi\left(a+(b+S)n^{-1/2}\right) \right\}}\\  
&  \xrightarrow[n \to \infty]{} &
\frac{ \pi(a)\textrm{E}^{*} \left\{I_{(-\infty, -k)}(S) \right\}}
{\pi(a)}  \\
&=& \Phi(-k),
\end{eqnarray*}
as claimed.
Specifically, n Figure \ref{fig:alpha_for_p}, we have $\pi(\theta)$ defined as the ``mixture'' prior (defined in equations \ref{eq:mixprior} and \ref{eq:norm_g0}) and set $a=0$ and $b=1.645$, with 4 different values of $k$: $k=1.645$, $k=1.282$, $k=0.842$, $k=0.524$, such that $\Phi(-1.645)=0.05$, $\Phi(-1.282)=0.10$, $\Phi(-0.842)=0.20$, and $\Phi(-0.524)=0.30$.
}

\end{document}